\documentclass[letterpaper,12pt]{article}
\usepackage{amsfonts,amsmath,amsthm,graphicx,enumerate,amssymb,amscd,youngtab,tikz,pgf,mathrsfs,pgfplots}
\usetikzlibrary{shapes.geometric, arrows}
\usepackage{pdfpages}
\pgfplotsset{compat=1.15}
\definecolor{uuuuuu}{rgb}{0,0,0}
\definecolor{qqqqff}{rgb}{0,0,1}
\definecolor{ududff}{rgb}{0.30196078431372547,0.30196078431372547,1}
\definecolor{qqttcc}{rgb}{0,0.2,0.8}
\definecolor{ffqqqq}{rgb}{1,0,0}
\definecolor{cqcqcq}{rgb}{0.7529411764705882,0.7529411764705882,0.7529411764705882}
\definecolor{amethyst}{rgb}{0.6, 0.4, 0.8}
\definecolor{amber}{rgb}{1.0, 0.49, 0.0}
 \definecolor{cygreen}{rgb}{0.0, 0.5, 0.0}
\usepackage{array,multirow,makecell}
\setcellgapes{1pt}
\makegapedcells
\newcolumntype{R}[1]{>{\raggedleft\arraybackslash }b{#1}}
\newcolumntype{L}[1]{>{\raggedright\arraybackslash }b{#1}}
\newcolumntype{C}[1]{>{\centering\arraybackslash }b{#1}}
\theoremstyle{definition}

\theoremstyle{remark}

\newcommand{\R}{\mathbb R}

\newcommand{\N}{\mathbb N}

\begin{document}

\title{The Cubic Equation Made Simple}

\author{Abdel Missa\footnote{Department of Finance, Jacksonville University. Email: amissa@ju.edu.}\footnote{MarketCipher Partners, Chief Investment Officer. Email: pm@market-cipher.com},  Chrif Youssfi \footnote{MarketCipher Partners, Quantitative Research. Email: cyoussfi@market-cipher.com}}

\maketitle

\newpage

\tableofcontents

\newpage
\begin{center}
\textbf{Abstract}\\
\end{center}
This article introduces an intuitive function MY that simplifies solving cubic equations without venturing into the complex space. To many, it’s quite strange that cubic root(s) are expressed using trigonometric functions in the three-real-roots case versus real-radicals in the one-real-root case. Yet, the MY function provides a different perspective to this oddity and shows that the transition between the two worlds is actually smooth.  \\
\\
Although $MY$’s behavior resembles power functions such as  $x^{\frac{1}{2}}$, $x^{\frac{1}{3}}$ and $x^{\frac{2}{5}}$, it cannot be expressed in real radicals. That said, we succeeded in proving an accurate closed form algebraic approximation of MY. So yes, casus irreducibilis still holds, but real radicals can get you very close!\\
\section{Introduction}
The cubic equation holds a special place in the history of mathematics.  In the early 16th century, the cubic formula was discovered independently by Niccolò Fontana Tartaglia and Scipione del Ferro. Italy at the time was famed for intense mathematical duels. In 1535, Tartaglia was challenged by Antonio Fior, del Ferro's student, with Tartaglia winning the famous contest. Gerolamo Cardano then persuaded Tartaglia to share his method with him, promising to not reveal it without giving Tartaglia time to publish. Once Cardano learned about Del Ferro's work, which predated Tartaglia's, he decided that his promise could be legitimately broken and published the method in his book “Ars Magna" in 1545. This led to a dramatic decade-long feud between Tartaglia and Cardano [1]. \\
\\
Later on, other methods were developed including a trigonometric solution for cubic equations with three real roots by François Viète (René Descartes expanded on Viète's work) [2]. Joseph Louis Lagrange followed with a new uniform method to solve lower degree (less than 5) polynomial equations including the cubic [3]. In 1683, Ehrenfried Walther von Tschirnhaus [4] proposed a new approach using the Tschirnhaus transformation. In addition, the authors identified other more recent methods and approaches to solving the cubic equation [5][6][7][8] and developed one of their own [9].\\
\\
In this article, a different approach to solving general cubic equations is discussed by introducing a new function MY that provides the real roots uniformly under various cases. MY is then expressed in closed form and hypergeometric form. We then proceed by using an algebraic iteration method that converges globally towards $MY$. \\
\\
While casus irreducibilis [10] for cubic equations states that real-valued roots of irreducible cubic
polynomials cannot be expressed in radicals without introducing complex numbers, we came really close with real radicals.  Finally, the article is concluded by discussing many of the unique properties of $MY$. 
\section {Canonical function}
Without loss of generality, let’s consider the depressed cubic equation\footnote{For a general cubic equations: $ax^3+bx^2+cx+d=0$ where $a \neq 0$, $b$, $c$ and $d$ are real numbers,  a change of variable $y=x+\frac{b}{3a}$ leads to the depressed.}:
\[y^3+p y +q=0 \qquad (1)\]
Where $p$ and $q$ are real numbers.\\
\\
Outside the trivial special cases of p=0 or q=0, equation (1) can be transformed to a canonical form:
\[\frac{z^3+ z^2}{2} = t\qquad (2) \]
Two transformations can be used:
\begin{enumerate}
\item \underline{Transformation 1}: A change of variable $z=\frac{q}{p y}$ leads to: $t=-\frac{q^2}{2p^3}$.
\item \underline{Transformation 2}: When $p<0$, a change of variable $z=\frac{y}{\sqrt{-3p}}-\frac{1}{3}$ leads to: \[ t=\frac{1}{27}-\frac{q}{2\sqrt{-27p^3}}\]
\end {enumerate} 
The canonical function $f$ is defined in $\R$ as: \[f: z \longmapsto \frac{z^3+ z^2}{2} \]
Using the sign of the derivative, there are three intervals where $f$ is monotonic:
\begin{enumerate}
\item $]-\infty,-\frac{2}{3}]$, $f$ is increasing from $-\infty$ to $\frac{2}{27}$. The point $M(-\frac{2}{3},\frac{2}{27})$ is a local maximum.
\item $[-\frac{2}{3},0]$, $f$ is decreasing from $\frac{2}{27}$ to $0$. The point $O(0,0)$ is a local minimum.
\item $[0,+\infty[$, $f$ is increasing from $0$ and $+\infty$.
\end{enumerate} 
Together, these properties define the number of roots of the canonical equation $f(z)= x$ (see Figure 1):
\begin{enumerate}
\item Scenario 1: $x>\frac{2}{27}$ , there is a unique solution that is above $\frac{1}{3}$.
\item Scenario 2: $x<0$, there is a unique negative solution.
\item Scenario 3: $0\leq x \leq \frac{2}{27}$ there are three real solutions (two of which may coincide). One root is positive and the other two are negative . 
\end{enumerate}
\begin{figure}[!h]
\centering
\begin{tikzpicture}
   \begin{axis}[ xmin=-1.7, xmax=1.7, ymin=-0.12, ymax=0.18,
     xscale=1.6, yscale=1.05]
     \draw[->, line width=1.1pt] (-1.5, 0) -- (1, 0) node[right] {$z$};
  	\draw[->, line width=1.1pt] (0, -0.18) -- (0, 0.14) node[above] {$y$};
	\draw (0.8,0.12) node[right, cygreen]{$ \textrm{Scenario 1}$};
	\draw (0.8,0.05) node[right, cygreen]{$ \textrm{Scenario 3}$};
	\draw (0.8,-0.08) node[right, cygreen]{$ \textrm{Scenario 2}$};
     \draw (0,0) node[below left]{$0$};
     \draw (0.33,2/27) node[above left, blue]{$f$};
     \addplot[blue, line width=0.8pt, samples=100, smooth, domain={-1.16}:{0.46}]plot (\x, { \x^3 + \x^2)*0.5 });
    	\addplot[black, dashed, samples=100, smooth, domain={-1.2}:{0.7}]plot (\x, {0.12 } );
     \addplot[black, dashed, samples=100, smooth, domain={-1.2}:{0.7}]plot (\x, { 0.05 } );
     \addplot[black, dashed, samples=100, smooth,  domain={-1.2}:{0.7}]plot (\x, { -0.08 } );
 	\draw (-2/3,2/27) node[above] {$M(-\frac{2}{3},\frac{2}{27})$};
  	\begin{scriptsize}
		\draw [fill=black] (0,0) circle (2pt);
		\draw [fill=black] (-2/3,2/27) circle (2pt);
	\end{scriptsize}
\end{axis}
\end{tikzpicture}
\caption{Canonical function $f$}
\label{Solution1}
\end{figure}
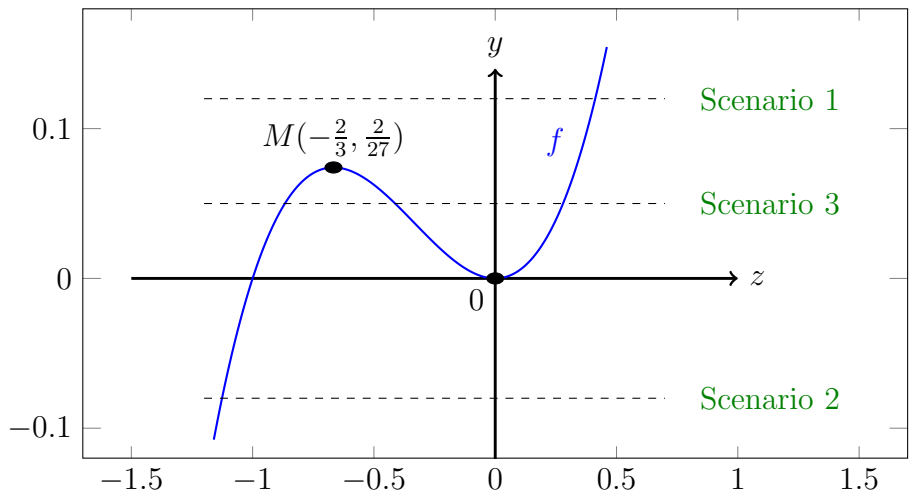
\section{Geometric intuition behind proposed method}
To provide the intuition behind our proposed method for solving cubic equations, notice that the inflection point $I\left(-\frac{1}{3}, \frac{1}{27}\right)$ is also a symmetry point. This can be expressed analytically as: 
\[f\left(-\frac{2}{3}-z\right)=\frac{2}{27}-f(z) \qquad \textrm{for all} \qquad z \in \R \]
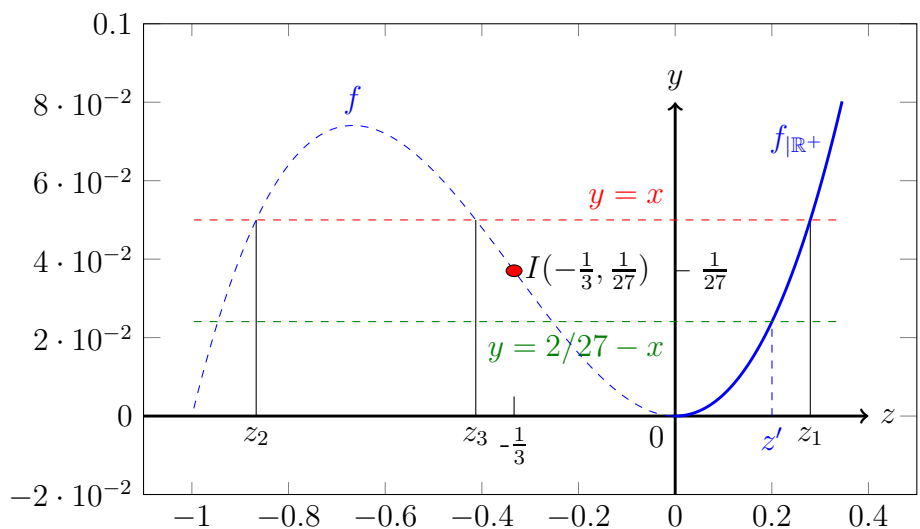
\begin{figure}[!h]
\center
\begin{tikzpicture}
   \begin{axis}[ xmin=-1.1, xmax=0.5, ymin=-0.02, ymax=0.1,
     xscale=1.5, yscale=1.1, restrict x to domain={-1}:0.35]
     	\draw[->, line width=1.1pt] (-1.1, 0) -- (0.4, 0) node[right] {$z$};
  	\draw[->, line width=1.1pt] (0, -0.1) -- (0, 0.08) node[above] {$y$};
     \draw (0.27955689, 0) node[below]{$z_1$ }--(0.27955689, 0.05);
     \draw (-0.866951318, 0) node[below]{$z_2$ }--(-0.866951318, 0.05);
     \draw (-0.412605572, 0) node[below]{$z_3$ }--(-0.412605572, 0.05);
	\draw [dashed, blue](0.200284651, 0) node[below]{$z'$ }--(0.200284651, 	2/27-0.05);
     \draw (0,0) node[anchor=north east]{$0$};
     \draw (0,1/27) --(0.03,1/27) node[right]{$\frac{1}{27}$};
     \draw (-1/3,0) node[below]{-$\frac{1}{3}$}; \draw (-1/3,0) -- (-1/3,0.005);
     \draw (-2/3,2/27) node[above, blue]{$f$};
     \draw (0.33,0.07) node[left, blue]{$f_{|\R^+}$};
     \draw (0,0.05) node[above left, red]{$y=x$};
     \draw (0,2/27-0.05) node[below left, cygreen]{$y=2/27-x$};
     \addplot[blue, samples=1000, smooth, line width=1.1pt, unbounded coords=discard, restrict x to domain={-0.01}:0.35]plot (\x, {(\x^3+ \x^2)*0.5) });
      \addplot[blue, dashed, samples=1000, smooth, unbounded coords=discard, restrict x to domain={-1}:0.35]plot (\x, { (\x^3+\x^2)*0.5 });
	\addplot[red, dashed, samples=1000, smooth, unbounded coords=discard, restrict x to domain={-1}:0.35]plot (\x, {0.05});
	\addplot[cygreen, dashed, samples=1000, smooth, unbounded coords=discard, restrict x to domain={-1}:0.35]plot (\x, {2/27-0.05});
	\draw (-1/3,1/27) node[right] {$I(-\frac{1}{3},\frac{1}{27})$};
     \begin{scriptsize}
		\draw [fill=red] (-1/3,1/27) circle (2pt);
	\end{scriptsize}
\end{axis}
\end{tikzpicture}
\caption{Solving of $f(z)=0.05$, three real roots}
\label{Solution2}
\end{figure}
Therefore, the problem of solving the equation $f(z)=x$ for $x \in \R$ is reduced to solving equations $f_{|\R^+}(z)=a$ for positive real numbers $a$. Let's illustrate this point geometrically by considering the following example $f(z)=x$ where $x=0.05$ (see Figure 2). This equation has three real solutions $z_1$, $z_2$ and $z_3$.
\begin{enumerate}
\item Construct the restricted curve $f_{|\R^+}$, as well as the lines $y=x$ and $y=\frac{2}{27}-x$. The abscissas of the intersections of $f_{|\R^+}$  with these two lines are respectively $z_1$, the positive root of the equation, and $z'$. 
\item  Let $z_2$ be the reflection of $z'$ with respect to $-\frac{1}{3}$: $z_2=-2/3-z'$. Using the symmetry property, $z_2$ is a root of the equation.
\item Let $z_3$ be the unique negative point such as the distance between $z_1$ and $-\frac{1}{3}$ is the same as the distance between $z'$ and $z_3$. In other words $z_3=-1/3-z'-z1$. Therefore $z_1+z_2+z_3=-1$. Using Vieta's formula, $z_3$ is the third root. 
\end{enumerate}

\section{$MY$ function definition}
The restriction, $f_{|\R^+}$, of $f$ to $\R^+$ is striclty increasing, continuous with  
$f_{|\R^+}(0)=0$ and $\lim\limits_{z \to +\infty} f_{|\R^+}(z)=+\infty$.
Therefore it is bijective from  $\R^+$ to  $\R^+$ and admits a reciprocal function:
\[MY: x \longmapsto MY(x)=f_{|\R^+}^{-1}(x) \qquad  x \in \R^+ \]
the graph of $MY$ is symmetrical to the graph of $f_{|\R^{+}}$ with respect to the line $y=x$ (See Figure 3). \\
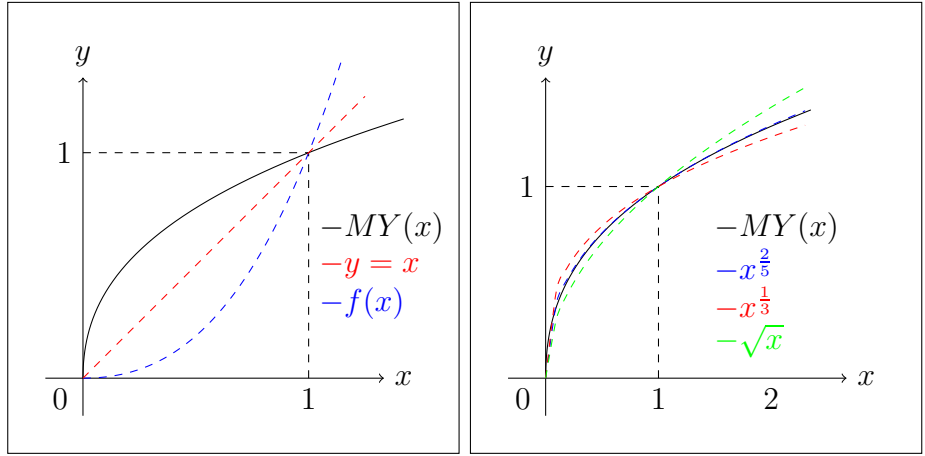
\begin{figure}[!h]
\centering
\begin{tikzpicture}
	\draw(-1,-1) rectangle (5,5);
     \draw[->] (-0.5, 0) -- (4, 0) node[right] {$x$};
  	\draw[->] (0, -0.5) -- (0, 4) node[above] {$y$};
     \draw [dashed](3, 0)node[below] {$1$} -- (3, 3);
	\draw [dashed](0, 3)node[left] {$1$} -- (3, 3);
     \draw [dashed](-0.3, 0)node[below] {$0$};
  	\draw[scale=3, domain=0:1.15,  dashed, smooth, variable=\x, blue] plot ({\x}, {(\x^3+\x^2)*0.5});
  	\draw[scale=3, domain=0:1.15, smooth, variable=\y, black]  plot ({(\y^3+\y^2)*0.5}, {\y});
  	\draw[scale=3, domain=0:1.25,  dashed, smooth, variable=\x, red]  plot ({\x}, {\x});
  	\draw (3, 2) node[right]{$\color{black} {- MY(x)}$};
  	\draw (3,1.5) node[right]{$\color{red}  {- y=x}$};
  	\draw (3,1) node[right]{$\color{blue}  {- f(x)}$};
  \end{tikzpicture}
\begin{tikzpicture}
	\draw(-1,-1) rectangle (5,5);
     \draw[->] (-0.5, 0) -- (4, 0) node[right] {$x$};
     \draw [dashed](1.5, 0)node[below] {$1$} -- (1.5, 2.55);
	\draw [dashed](0, 2.55)node[left] {$1$} -- (1.5, 2.55);
     \draw [dashed](-0.3, 0)node[below] {$0$};
     \draw (3, 0)node[below] {$2$};
  	\draw[->] (0, -0.5) -- (0, 4) node[above] {$y$};
  	\draw[scale=1.5, domain=0:1.4, smooth, variable=\x, black]  plot ({(\x^3+\x^2)*0.5}, {1.7*\x});
  	\draw[scale=1.5, domain=0:2.3,  dashed, smooth, variable=\x, blue] plot ({\x}, {1.7*\x^0.4});
  	\draw[scale=1.5, domain=0:2.3, dashed, smooth, variable=\x, red]  plot ({\x}, {1.7*\x^0.333});
  	\draw[scale=1.5, domain=0:2.3, dashed, smooth, variable=\x, green]  plot ({\x}, {1.7*\x^0.5});
  	\draw (2.1,2) node[right]{$\color{black} {- MY(x)}$};
  	\draw (2.1,1.5) node[right]{$\color{blue} {- x^{\frac{2}{5}}}$};
  	\draw (2.1,1) node[right]{$\color{red} {- x^{\frac{1}{3}}}$};
	\draw (2.1,0.5) node[right]{$\color{green} {- \sqrt{x}}$};
\end{tikzpicture}
\caption{{\bf Left}: $MY$ inverse of $f$ \qquad {\bf Right:} $MY$ versus power functions }
\end{figure}
\\
$MY$ is continuous, strictly increasing and infinitely differentiable. Its behavior resembles  power functions (See Figure 3): When $x$ is close to 0, $MY(x)$ is equivalent to $\sqrt{2x}$. For large values of $x$, $MY(x)$ behaves like $\sqrt[3]{2x}$ (see section 8, Properties of $MY$). $x^{\frac{2}{5}}$ is an upper bound for $MY$.\\
\\
Thanks to Cardano and Vieta's trignometric formulas  for cubic roots, $MY$ can be expressed in closed form. Define: 
\[u=x-\frac {1}{27}\]

\begin{enumerate}
\item For $x \in [\frac{2}{27},+\infty [$\\
\[MY\left(x\right)=-\frac{1}{3} +\sqrt[3]{u+\sqrt{u^2-\left (\frac {1}{27}\right)^2}}+\sqrt[3]{u-\sqrt{u^2-\left (\frac {1}{27}\right)^2}}\]
\item For $x \in [0, \frac{2}{27}[$\\
\[MY\left(x\right)=-\frac {1}{3}+\frac{2}{3}\cos \left (\frac{\arccos (27 u)}{3}\right)\]  
\end{enumerate}
Notice that MY is continuous and smooth at $x= \frac{2}{27}$ despite two totally different closed form expressions!
\section{Cubic roots expressed in $MY$}

\subsection{Solving the canonical equation}
Roots of the canonical equation $f(z)=x$ can be expressed in a simple form using $MY$:
\begin{enumerate}
\item For $x > \frac{2}{27}$, there is a unique real solution $z_1=MY\left(x\right)$. 
\item For $x <  0$, there is a unique real solution $z_1$.  Using the symmetry property:
\[f\left( -\frac{2}{3}-z_1\right)=\frac{2}{27}-x \qquad \textrm{it follows that:} \qquad
z_1=-\frac{2}{3}-MY\left(\frac{2}{27}-x\right) \]
\item For $ 0 \leq x \leq \frac{2}{27}$ there are three real solutions. $z_1=MY\left(x\right)$ is the positive solution.  The two other negative solutions are obtained, first, by  symmetry:
\[z_2=-\frac{2}{3}-MY\left(\frac{2}{27}-x\right) \]
second, by using Vieta's formula\footnote{Also, for $x \neq 0$ $z_3=2x/(z_1z_2)$} $z_1+z_2+z_3=-1$ : 
\[z_3=MY\left(\frac{2}{27}-x\right)-MY\left(x\right)-\frac{1}{3} \]
Given the sign of the three roots: $z_2, z_3 \leq 0 \leq z_1$. In addition:
\[z_3-z_2=\frac{1}{3}+2 MY\left(\frac{2}{27}-x\right)-MY\left(x\right)\]
Since $MY(x) \leq MY(\frac{2}{27})=\frac{1}{3}$:
\[z_2 \leq z_3 \leq z_1 \qquad (4)\] 
\end{enumerate}

\subsection{Solving the depressed cubic}
Using the results from section 2 and section 5.1, we can express the roots of the depressed equation:
\[y^3+p y +q=0 \qquad (1)\]
Indeed, equation (1) can be transformed to a canonical form using either Transformation 1 or Transformation 2. Interestingly, each transformation leads to a different expression of the roots. 
In particular, when $p<0$ both transformations can be applied which leads to useful equalities. For this purpose, let’s define:
\[\xi=\frac{3 q}{2 p}\sqrt{-\frac{3}{p}}\]
We distinguish 4 cases:\\
\\
{\bf Case 1: $p =0$}. There is a unique real solution $\alpha =-\sqrt[3]{q}$.\\
\\
{\bf Case  2: $p > 0$}. There is a unique real solution:
\[\alpha =\frac{q}{p \left(-\frac{2}{3}-MY\left(\frac{2}{27}+\frac{q^2}{2p^3}\right)\right)}\]
{\bf Case  3: $p < 0$ and $|\xi| >1$}. There is a unique real solution with two expressions:
\[\alpha =-\sqrt{-\frac{p}{3}}\left(3 MY\left(\frac{1+|\xi|}{27}\right) +1\right)=\frac{q}{p MY\left(-\frac{q^2}{2p^3}\right)}\]
{\bf Case  4: $p < 0$ and $-1\leq \xi \leq 1$}. There are three real solutions that can be expressed using Tranformation 2 as:
\[\alpha=\sqrt{-\frac{p}{3}}\left(3 MY\left(\frac{1+\xi}{27}\right) +1\right)\]
\[\beta=3\sqrt{-\frac{p}{3}}\left(MY\left(\frac{1-\xi}{27}\right) -MY\left(\frac{1+\xi}{27}\right) \right)\]
\[\gamma=-\sqrt{-\frac{p}{3}}\left(3 MY\left(\frac{1-\xi}{27}\right) +1\right)\]
As a result of (4):
\[\gamma \leq \beta \leq \alpha\]
When $q \neq 0$, these roots can be expressed differently using Transformation 1:
\[\alpha'=\frac{q}{p MY\left(-\frac{q^2}{2p^3}\right)}\]
\[\beta'=\frac{q}{p \left(MY\left(\frac{2}{27}+\frac{q^2}{2p^3}\right)-MY\left(-\frac{q^2}{2p^3}\right)-\frac{1}{3}\right)}\]
\[\gamma'=\frac{-q}{p \left(\frac{2}{3}+MY\left(\frac{2}{27}+\frac{q^2}{2p^3}\right)\right)}\]
Alternatively, $\beta'$ could be derived from any of Vieta's formulas for equation $(1)$. For example: \[\alpha'+\beta'+\gamma'=0\]
In addition, using the order of the three roots, when $q < 0$:
\[\alpha'=\alpha \qquad \beta'=\gamma \qquad \gamma'=\beta \]
and when $q>0$: 
\[\alpha'=\gamma \qquad \beta'=\alpha \qquad \gamma'=\beta \]
Naturally, the roots coincide with the trigonemetric roots provided by François Viète:
\[t_k=2\sqrt{-\frac{p}{3}}\cos\left(\theta_k\right) \qquad \textrm{where}\qquad \theta_k= \frac{\arccos(\xi)}{3}-\frac{2 k\pi}{3}\qquad \textrm{for } k=0,1,2\]
Since $\theta_0 \in [0, \pi/3]$, $\theta_1 \in [-2\pi/3, -\pi/3]$:  $\theta_2 \in [-4\pi/3, -\pi]$:  
\[t_2 \leq t_1 \leq t_0\]
Therefore
 \[\alpha=t_0 \qquad \beta=t_1 \qquad \gamma=t_2\]
\section{Approximation of $MY$ using radicals}
\subsection{Casus irreducibilis}
One of the oddities of solving general cubic equations using radicals is the absolute  requirement to use complex numbers in the irreducible case. In other words, despite the three roots being real for irreducible cubic polynomials, they cannot be expressed as radicals of real numbers [10]. We note of course the alternate solution provided by Viète, which bypasses the use of complex  numbers by introducing trigonometric functions. This is reflected indirectly in the provided closed form of the function $MY$, which is a hybrid of an algebraic function for $x \geq \frac{2}{27}$ and a trigonometric function for $x<\frac{2}{27}$.\\
\\
In the following section, we provide an accurate closed form algebraic approximation of $MY$. In doing so, we offer a new perspective on solving the cubic equation including in the irreducible case. 
 
\subsection{Fixed point iteration}
Assume that $z=MY(x)$. z satisfies the equation: 
\[z^3+z^2=2x \qquad (5)\] 
This equation can be exploited in two ways. First, by factorizing $(5)$ we obtain:
\[z=\sqrt{\frac{2x}{1+z}} \qquad (6)\]
Second, we can complete the cubic: 
\[\left(z+\frac{1}{3}\right)^3=2x+\frac{1}{27}+\frac{z}{3} \qquad (7)\]
When substituting $z$ in the right hand side of $(7)$ by using (6):
\[z=G(x,z) \qquad (8)\]
Where the function $G$ is defined by:
\[ G\left(x,y\right)=\sqrt[3]{2 x+ \frac{1}{27}+\frac{1}{3}\sqrt{\frac{2 x}{1+y}}}-\frac{1}{3}\qquad (9)\]
Therefore z is a fixed point of $G(x,.)$, inspiring a fixed point iteration.\\
\\
The following theorem proves the convergence of this method towards $MY$. More importantly, it provides right away an accurate closed form algebraic approximation of $MY$ across $\R^+$.\\
\\
{\bf Convergence theorem:}\\
\\
Define the sequence of functions $(M_n(.))_{n \in \N}$ defined by:
\[M_0(x)=G(x,x^{\frac{2}{5}}) \qquad  \textrm{and} \qquad M_{n+1}(x)=G(x,M_n(x))\]
\begin{enumerate}
\item For all positive real numbers $x$:
\[|M_{0}(x)-MY(x)|<C_0 \qquad \textrm{where } \qquad C_0=1.4408 \: 10^{-3}\]
\item For all strictly positive real numbers $x$:
\[\left|\frac{M_{0}(x)}{MY(x)}-1\right|<{C_0}' \qquad \textrm{where } \qquad {C_0}'=1.1527 \: 10^{-2}\]
\item For all positive real numbers $x$:
\[|M_{n}(x)-MY(x)|<\frac{C_0}{K^n} \qquad \textrm{where } \qquad  K=25.05\]
\item This sequence converge uniformly to MY over $\R^+$.
\\
\end {enumerate}
\newpage
In order to demonstrate the convergence theorem, we start with the following lemma:\\
\\
{\bf Lemma:}
\begin{enumerate}
\item For all positive real numbers $x$ and $y$:\[\left|\frac{\partial G}{\partial y}\left(x,y\right) \right| < C_1 \qquad \textrm{where}\qquad  C_1 \approx \frac{1}{21.2398}\]
\item For all positive real numbers $x$ : \[\left|\frac{\partial G}{\partial y}\left(x,MY\left(x\right)\right)\right | < C_2 \qquad \textrm{where} \qquad C_2 \approx \frac{1}{30.5475}\]
\end {enumerate}
{\bf Proof of lemma:}\\
\\
$G(x,.)$ is decreasing with respect to $y$ and $G(.,y)$ is increasing with respect to $x$. Define: 
\[ z=MY(x) \qquad (10)\] 
z satisfies:
\[ G\left(x,z\right)=z \qquad (11)\]
For any $x>0$, $G(x,.)$ is derivable and for $y \geq 0$ : \[ \frac{\partial G}{\partial y}\left(x,y\right)=-\frac{\sqrt{2x}}{18}\left( \left(2 x+ \frac{1}{27}\right) (1+y)^{\frac{9}{4}}+\frac{\sqrt{2x}}{3}(1+y)^{\frac{7}{4}}\right)^{-\frac{2}{3}} \qquad (12)\]
\begin{enumerate}
\item Upper bound for $\left|\frac{\partial G}{\partial y}\left(x,.\right)\right|$\\
$\left|\frac{\partial G}{\partial y}\left(x,.\right)\right|$ is strictly decreasing and convex. Therefore its maximum $C_1$ is reached at $y=0$.\\
\\
Define $t=\sqrt{2x}$. $C_1$ is also the maximum of:
\[h(t)=\frac{t}{18}\left( t^2+ \frac{1}{27}+\frac{t}{3}\right)^{-\frac{2}{3}}=\frac{1}{18}\left( t^{\frac{1}{2}}+ \frac{t^{-\frac{3}{2}}}{27}+\frac{t^{-\frac{1}{2}}}{3}\right)^{-\frac{2}{3}}\]
Let $v=t^{-\frac{1}{2}}$:
\[h(t)=g(v)=\frac{1}{18}\left( \frac{1}{v}+ \frac{v^3}{27}+\frac{v}{3}\right)^{-\frac{2}{3}}\]
$C_1$, is obtained by setting the derivative of $g$ to $0$, reached at $v_0$:
\[v_0=\sqrt{\frac{-1+\sqrt{5}}{2}} \]
and
 \[ C_1 \approx \frac{1}{21.2398}\]
\item Upper bound for $\left|\frac{\partial G}{\partial y}\left(x,MY(x)\right)\right|$\\
\\
For $x>0$ define $z=MY(x)$:
\[\left| \frac{\partial G}{\partial y}\left(x,z\right)\right| =\frac{\sqrt{2x}}{18(1+z)^{\frac{3}{2}}}\left( \left(2 x+ \frac{1}{27}\right) +\frac{1}{3}\sqrt{\frac{2x}{1+z}}\right)^{-\frac{2}{3}}  \]
Recall that:
\[z=\sqrt{\frac{2x}{1+z}} \qquad \textrm{and} \qquad (z+\frac{1}{3})^3=2x+\frac{1}{27}+\frac{z}{3}\] 
Therefore
\[ \left|\frac{\partial G}{\partial y}\left(x,z\right)\right| =\frac{z}{18 (1+z) \left(z+\frac{1}{3}\right)^2} \qquad (13)\]
Or:
\[ \left|\frac{\partial G}{\partial y}\left(x,z\right)\right|= \frac{1}{18 }\left(z^2+\frac{5}{3}z + \frac{7}{9}+\frac{1}{9z}\right)^{-1}\]
Setting the derivative to $0$ leads to the maximum $C_2$, reached at $z$ such that:
\[2 z+\frac{5}{3}-\frac{1}{9z^2}=0\]-
Notice that:
\[2 z+\frac{5}{3}-\frac{1}{9z^2}=\frac{(3z+1)(6 z^2+3 z-1)}{9 z^2}\]
The only positive solution is: \[z=\frac{-3+\sqrt{33}}{12}\] leading to: \[C_2 \approx \frac{1}{30.5475}\]
\end{enumerate}
{\bf Proof of Convergence theorem}
\begin{enumerate}
\item For $x>0$, define $z=MY(x)$. \\
Since $M_0(x)=G(x,x^{\frac{2}{5}})$ and $z=G(x,z)$:
\[M_0(x)-MY(x)|=\left|G(x,x^{\frac{2}{5}})-z \right|\] 
As proven below in the Properties of $MY$ section $z \leq x^{\frac{2}{5}}$. \\
Since $\left|\frac{\partial G}{\partial y}(x,.)\right| $ is decreasing with respect to $y$ over the interval $[z,x^{\frac{2}{5}}]$:
\[ \left|M_0(x)-MY(x)\right| = \left|G(x,x^{\frac{2}{5}})-G(x,z) \right| \leq U(z)=\left|\frac{\partial G}{\partial y}(x,z)\right|\left|(x^{\frac{2}{5}}-z)\right| \qquad (14)\]
Recall:
\[x=(z^3+z^2)/2 \qquad \textrm{and } \qquad \left|\frac{\partial G}{\partial y}(x,z)\right|=\frac{z}{18(1+z)\left(z+\frac{1}{3}\right)^{2}}\]
Therefore
\[U(z)= \frac{z \left(\left( \frac{z^3+z^2}{2}\right)^{\frac{2}{5}} -z\right)}{18(1+z)\left(z+\frac{1}{3}\right)^{2}}\]
Write \[U(z)=A.B \qquad (15)\] where:
\[A=\frac{z^2(1+\frac{1}{z})}{18\left(z+\frac{1}{3}\right)^{2}}=\frac{1}{18}\left(1+\frac{1}{3}\frac{\left(z-\frac{1}{3}\right)}{\left(z+\frac{1}{3}\right)^{2}}\right)\]
\[B=\frac{\left( \frac{z^{1/2}+z^{-1/2}}{2}\right)^{\frac{2}{5}} -1}{(1+z)(1+\frac{1}{z})} \]
First, when $z\leq \frac{1}{3}$ $A\leq \frac{1}{18}$.  If $z >\frac{1}{3}$, define  $y=z-\frac{1}{3}>0$.
\[A=\frac{1}{18}\left(1+\left(y^{\frac{1}{2}}+\frac{2}{3}y^{\frac{-1}{2}}\right)^{-2}\right)\]
Using the derivative, A is maximal when $y=\frac{2}{3}$ or $z=1$. Therefore
\[A \leq \frac{1}{16} \qquad (16)\] 
Second, define:
\[w=\frac{z^{1/2}+z^{-1/2}}{2} \geq 1\]
Notice that 
\[(1+z)(1+\frac{1}{z})=4w^2 \qquad \textrm{and} \qquad B=\frac{w^\frac{2}{9}-1}{4w^2}\]
When $w=1$ $B=0$. If $w>1$ define $\xi=w^\frac{2}{9}-1$. A simple algebra leads to:
\[B=\frac{1}{4} \left(\xi^\frac{-2}{9}+\xi^\frac{7}{9}\right)^{-\frac{9}{2}}\]
Using the derivative $B$ is maximal when $\xi=\frac{2}{7}$. Therefore:
\[B \leq \frac{1}{43.37886} \qquad (17)\]
As a result:
\[ U(z) \leq \frac{1}{694.061782} \approx  1.4408 \ 10^{-3} \]
Going back to $(14)$, for all positive real numbers $x$:
\[\left| M_0(x)-MY(x) \right|<C_0 \qquad (18)\] 
Where
\[C_0=\frac{1}{694.061782}\approx  1.4408 \: 10^{-3}\]
\item Keeping the notation $z=MY(x)$ and using (14) and (15):
\[\left| \frac{M_0(x)}{MY(x)}-1 \right|<\frac{U(z)}{z}={\tilde A}(z) B \]
Where
\[{\tilde A}(z)=\frac{A}{z}=\frac{z+1}{18\left(z+\frac{1}{3}\right)^{2}}\]
The derivative ${\tilde A}'$ of ${\tilde A}$ is
\[{\tilde A}'(z)=-\frac{\left(z+\frac{5}{3}\right)}{18\left(z+\frac{1}{3}\right)^{3}}<0 \qquad \textrm{for }z\geq 0\]
Therefore ${\tilde A}$ is decreasing on $\R^+$ and reaches its maximum at $z=0$. Then
 \[{\tilde A}(z) \leq {\tilde A}(0)=\frac{1}{2}\] 
This combined with (17) leads to:
\[\left| \frac{M_0(x)}{MY(x)}-1 \right|<{\tilde A}(z) B \leq {C_0}'= 1.1527 \: 10^{-2} \]

\item For $n \in \N$  and any positive real number $x$, define $z=MY(x)$ and $y=M_n(x)$. Since $M_{n+1}(x)=G(x,y)$ and $z=G(x,z)$ :
\[\left|M_{n+1}(x)-z\right|= \int_{y}^{z}\left|\frac{\partial G}{\partial y}\left(x,u\right)\right|du \qquad (19)\]
Since $\left|\frac{\partial G}{\partial y}\left(x,.\right)\right|$ is positive and convex, the integral in $(19)$ is lower than the area of the trapezoid:
\[ \left|M_{n+1}(x)-z\right|= \leq \frac{1}{2}\left|y-z\right|)\left(\left|\frac{\partial G}{\partial y}\left(x,y\right)\right |+\left|\frac{\partial G}{\partial y}\left(x,z\right)\right |\right)\]
Using results from the lemma: 
\[ \left|M_{n+1}(x)-z\right| \leq \left|M_{n}(x)-z\right|\frac{C_1+C_2}{2} \qquad (20)\]
Or
\[ \left|M_{n+1}(x)-z\right| \leq \frac{\left|M_{n}(x)-z\right|}{K} \qquad (21)\]
Where
\[ K=\frac{2}{C_1+C_2}\approx 25.0572 \]
It follows from $(A.8)$ and $(A.11)$  that: 
\[ \left|M_{n}(x)-MY(x)\right| \leq \frac{C_0}{K^n} \qquad (22)\]
\item As a  natural consequence of  $(22)$, the sequence $(M_n)_{n \in \N}$ converges uniformly to MY. 
\end{enumerate}
\subsection{Examples:}
\subsubsection {Numerical examples for the approximation of MY}
To evaluate the efficiency of the algorithm described above, we consider two examples:\\	
\\	
{\bf Example 1:} $x=0.01$, $MY(x)$ closed form approximate value to 10th decimal digit is $0.1328694292$.\\
\\
\begin{tabular}{|R{1.5cm}|R{3.5cm}|R{3.5cm}|R{3cm}|}
\hline Iteration & $M_n(x)$ & $\left|M_n(x)-MY(x)\right|$ & $\left|\frac{M_n(x)}{MY(x)}-1\right|$   \\
\hline 0 & 0.1321129198 & 7.57  $10^{-04}$  & 5.69  $10^{-03}$\\
\hline 1 & 0.1328921191 & 2.27  $10^{-05}$  & 1.71  $10^{-04}$ \\
\hline 2 & 0.1328687489 & 6.80  $10^{-07}$  & 5.12  $10^{-06}$\\
\hline 3 & 0.1328694495 & 2.04  $10^{-08}$  & 1.53  $10^{-07}$\\
\hline 4 & 0.1328694285 & 6.11  $10^{-10}$  & 4.60  $10^{-09}$\\
\hline 5 & 0.1328694292 & 1.83  $10^{-11}$  & 1.38  $10^{-10}$\\
\hline 
\end{tabular}
\\	
\\
\\
{\bf Example 2:} $x=1000$,  $MY(x)$ closed form approximate value to 10th decimal digit is  $12.2745406200$. \\
\\
\begin{tabular}{|R{1.5cm}|R{3.5cm}|R{3.5cm}|R{3cm}|}
\hline Iteration & $M_n(x)$ & $\left|M_n(x)-MY(x)\right|$ & $\left|\frac{M_n(x)}{MY(x)}-1\right|$   \\
\hline 0 & 12.2735762826  & 9.64  $10^{-04}$  & 7.86  $10^{-05}$\\
\hline 1 & 12.2745409317  & 3.12  $10^{-07}$  & 2.54  $10^{-08}$ \\
\hline 2 & 12.2745406200  & 6.85  $10^{-11}$  & 5.58  $10^{-12}$ \\
\hline 
\end{tabular}
\\
\\
\\
Note that when $x<2/27$ (example 1), the closed form expression of $MY$ is trigonometric, and when $x \geq \frac{2}{27}$ (example 2) , the closed form expression is in radicals. Yet, the global algebraic (in radicals) approximation provided works well for both scenarios. 
\subsubsection {Numerical examples for solving cubic equations}

{\bf Example 1:} $x^3+x+1=0$.\\
Here $p=1$ and $q=1$. The results of section 5.2, case 2 apply. There is a unique real solution $\alpha$. Using $MY$ closed form:
\[ \alpha \approx-0.6823278038 \] 
$\alpha_n $ is the estimate of $\alpha$ using $ n $ iterations of the fixed point algorithm of  $MY$.\\
\\
\begin{tabular}{|R{0.5cm}|R{4cm}|R{3cm}|R{3cm} |}
\hline $n$ & $ \alpha_n $ & $ | \alpha_n- \alpha | $ &  $ \left|\frac{\alpha_n}{\alpha}-1 \right|$  \\
\hline 0 & -0.6823458163	& 1.80 $10^{-05}$	& 2.64 $10^{-05}$\\
\hline 1 & -0.6823274572	& 3.47 $10^{-07}$	& 5.08 $10^{-07}$\\
\hline 2 & -0.6823278105	& 6.67 $10^{-09}$	& 9.78 $10^{-09}$\\
\hline 3 & -0.6823278037	& 1.28 $10^{-10}$	& 1.88 $10^{-10}$\\
\hline 
\end{tabular}
\\
\\
\\
{\bf Example 2:} $x^3-3x+1=0$.\\
\\
Here $p=-3$ and $q=1$ which leads to $\xi =q\sqrt{-\frac{27}{4p^3}}= \frac{1}{2}$. The results of section 5.2, case 4 apply. There are three real roots $\alpha$, $\beta$ and $\gamma$. Using the closed form of $MY$:
\[ \alpha \approx 1.5320888862 \qquad \beta \approx 0.3472963553 \qquad \gamma\approx -1.8793852416 \] 
$\alpha_n $, $\beta_n$ and $\gamma_n$ are the respective estimates of $\alpha$, $\beta$ and $\gamma$  using $ n $ iterations of the fixed point algorithm of  $MY$.\\
\\
\begin{tabular}{|R{0.5cm}|R{4cm}|R{3cm}|R{3cm} |}
\hline $n$ & $ \alpha_n $ & $ | \alpha_n- \alpha | $ &  $ \left|\frac{\alpha_n}{\alpha}-1 \right|$  \\
\hline 0 & 1.5296764368	& 2.41 $10^{-03}$	& 1.57 $10^{-03}$\\
\hline 1 & 1.5321663348	& 7.74 $10^{-05}$	& 5.06 $10^{-05}$\\
\hline 2 & 1.5320864010	& 2.49 $10^{-06}$	& 1.62 $10^{-06}$\\
\hline 3 & 1.5320889660	& 7.79 $10^{-08}$	& 5.21 $10^{-08}$\\
\hline 4 & 1.5320888837	& 2.56 $10^{-09}$	& 1.67 $10^{-09}$\\
\hline 5 & 1.5320888863	& 8.21 $10^{-11}$	& 5.36 $10^{-11}$\\
\hline 
\end{tabular}
\\
\\
\begin{tabular}{|R{0.5cm}|R{4cm}|R{3cm}|R{3cm} |}
\hline $n$ & $ \beta_n $ & $ | \beta_n- \beta | $ &  $ \left|\frac{\beta_n}{\beta}-1 \right|$  \\
\hline 0 & 0.3476559549	& 3.60 $10^{-04}$	& 1.04 $10^{-03}$\\
\hline 1 & 0.3472848043	& 1.16 $10^{-05}$	& 3.33 $10^{-05}$\\
\hline 2 & 0.3472967260	& 3.71 $10^{-07}$	& 1.07 $10^{-06}$\\
\hline 3 & 0.3472963434	& 1.19 $10^{-08}$	& 3.42 $10^{-08}$\\
\hline 4 & 0.3472963557	& 3.82 $10^{-10}$	& 1.10 $10^{-09}$\\
\hline 5 & 0.3472963553	& 1.22 $10^{-11}$	& 3.53 $10^{-11}$\\
\hline 
\end{tabular}
\\
\\
\\
\begin{tabular}{|R{0.5cm}|R{4cm}|R{3cm}|R{3cm} |}
\hline $n$ & $ \gamma_n $ & $ | \gamma_n- \gamma | $ &  $ \left|\frac{\gamma_n}{\gamma}-1 \right|$  \\
\hline 0 & -1.8773323917	& 2.05 $10^{-03}$	& 1.09 $10^{-03}$\\
\hline 1 & -1.8794511391	& 6.59 $10^{-05}$	& 3.51 $10^{-05}$\\
\hline 2 & -1.8793831270	& 2.11 $10^{-06}$	& 1.13 $10^{-06}$\\
\hline 3 & -1.8793853094	& 6.79 $10^{-08}$	& 3.61 $10^{-08}$\\
\hline 4 & -1.8793852394	& 2.18 $10^{-09}$	& 1.16 $10^{-09}$\\
\hline 5 & -1.8793852416	& 6.99 $10^{-11}$	& 3.72 $10^{-11}$\\
\hline 
\end{tabular}
\\
\\
\\
Obviously, the goal here is not to provide the most optimal root-finding algorithm (see Newton  method, Secant method, Steffensen method,  Halley method, Laguerre method, Aberth-Ehrlich method, Durand-Kerner method, etc.). Instead, our aim is to merely shed light on an algorithm that approximates with real-radicals the function $MY$ and the roots especially in casus irreducibilis. 

\section{Hypergeometric representation}
The objective of this section is to express $MY$ using hypergeometric functions. Recall that $MY(x)$ is the unique positive solution of the equation:
\[z^3+z^2-2x=0\]
For $x>0$ we consider $y=1/z$, this implies that:
\[y^3-\frac{1}{6x}y-\frac{1}{4x}=0\]
We use the method provided by Zucker in [8]. Define:
\[p'=-\frac{1}{2x} \qquad q'=-\frac{1}{2x} \qquad \textrm{and} \qquad \Delta'=q'^2+p'^3\]
The Cardano formula expresses the root as $y=u+v$ with:
\[u=\left(-q'+\sqrt{\Delta'}\right)^{-\frac{1}{3}} \qquad v=\left(-q'-\sqrt{\Delta'}\right)^{-\frac{1}{3}}\]
Or:
\[u+v=(-q')^{-\frac{1}{3}}\left(\left(1+\sqrt{\frac{\Delta'}{q'^2}}\right)^{-\frac{1}{3}}+\left(1-\sqrt{\frac{\Delta'}{q'^2}}\right)^{-\frac{1}{3}}\right)\]
Using the identity:
\[(1+z)^{-2a}+(1-z)^{-2a}=2 F(a,a+\frac{1}{2},\frac{1}{2},z^2) \qquad \textrm{for }\qquad a=-\frac{1}{6}\]
where $F$ is the Gaussian hypergeometric function (analytically continued), we obtain:
\[u+v=(-q')^{-\frac{1}{3}}F\left(-\frac{1}{6}, \frac{1}{3}, \frac{1}{2},\frac{\Delta'}{q'^2}\right)\]
Using Kummer's transformation [11]: $F(a,b,c,z)=(1-z)^{-b}F(c-a,b,c,\frac{z}{z-1})$:
\[u+v=\left(\frac{2q'}{p'}\right)^{-\frac{1}{3}}F\left( \frac{1}{3},\frac{2}{3}, \frac{1}{2},\frac{\Delta'}{p'^3}\right)\]
Or:
\[u+v=3 F\left(\frac{2}{3}, \frac{1}{3}, \frac{1}{2},1-\frac{27 x}{2}\right) \qquad (23)\]
Using Kummer's transformation a second time:
\[u+v=3 \left(\frac{27 x}{2}\right)^{-\frac{2}{3}}F\left(\frac{1}{6}, \frac{2}{3}, \frac{1}{2},1-\frac{2}{27x}\right) \qquad (24)\]
If $x\leq \frac{2}{27}$, using formula $(23)$, $u+v$ is positive. Likewise when $x\geq \frac{2}{27}$, formula $(24)$ shows that $u+v$ is positive.\\
\\
Therefore
\[{\bf MY(x)=\frac{1}{3 F\left(\frac{1}{3},\frac{2}{3},  \frac{1}{2},1-\frac{27 x}{2}\right)}} \]
Notice that  
\[\lim\limits_{z \to 1}F\left(\frac{1}{3},\frac{2}{3},  \frac{1}{2},z\right)=+\infty \]
which is coherent with $MY(0)=0$.

\section{Properties of $MY$}
In Annex we prove all of the following properties:\\
\\
{\bf Equalities:}
\begin{enumerate}
\item For $x \in [\frac{2}{27},+\infty [$:
\[MY\left(x\right)= \frac{\sqrt[3]{2x\left(x-\sqrt{x\left (x-\frac {2}{27}\right)}\right)}}{\frac{1}{3}+\sqrt[3]{\left (x-\frac {1}{27} \right)-\sqrt{x\left (x-\frac {2}{27}\right)}}}\]
\item For all positive real numbers $x$, $MY$ satisfies the identity :
\[MY\left(x\right)\left( 3 MY\left(\sqrt{\frac{x}{54}}+\frac{1}{27}\right)+1\right)=\sqrt{6 x}\]
\item For all positive real numbers $x \neq 0 $:
\[MY(x)=\frac{1}{MY\left(\frac{x}{(MY(x))^5}\right) }   \]
\item For $0 \leq x \leq \frac{2}{27}$, the three roots of the equation $z^3+z^2=2x$ are:
\[z_1=MY(x)\]
\[z_2=-\frac{2}{3}-MY\left(\frac{2}{27}-x\right) =-\frac{1+MY(x)}{2}\left(1+\sqrt{\frac{1-3 MY(x)}{1+MY(x)}}\right)\]
\[z_3=MY\left(\frac{2}{27}-x\right)-MY(x)-\frac{1}{3} =-\frac{1+MY(x)}{2}\left(1-\sqrt{\frac{1-3 MY(x)}{1+MY(x)}}\right)\]
Consequently
\[MY\left(\frac{2}{27}-x\right) =\frac{1+MY(x)}{2}\left(1+\sqrt{\frac{1-3 MY(x)}{1+MY(x)}}\right)-\frac{2}{3}\]
\item \[\lim\limits_{x \to 0}\frac{MY\left(x\right)}{\sqrt{2x}}=1  \qquad \textrm{and} \qquad \lim\limits_{x \to +\infty}\frac{MY\left(x\right)}{\sqrt[3]{2x}}=1\]\\
Consequently for any $x \geq 0$
\[  \sqrt{x}=\lim\limits_{\epsilon \to 0^+} \frac{1}{\epsilon} MY\left(\frac{x \epsilon^2}{2}\right) \qquad \textrm{and} \qquad  \sqrt[3]{x}=\lim\limits_{\epsilon \to 0^+} \epsilon MY\left(\frac{x}{2\epsilon^3}\right)\]
\end{enumerate}

{\bf Inequalities:}
\begin{enumerate}
\item For all positive real numbers $x$:
\[\sqrt{\frac{2 x}{1+x^{\frac{2}{5}}}} \leq MY(x) \leq x^{\frac{2}{5}} \]
\item For all positive real numbers $x\neq 0$:
\begin{enumerate}
	\item If $a$ is a real number such that $0 \leq a \leq1$:
	\[MY\left(x^a\right) \geq \left(MY(x)\right)^a\]
	\item If $a$ is a real number such that $a \leq 0$ or $a \geq 1$:
	\[MY\left(x^a\right) \leq \left(MY(x)\right)^a\]
\end{enumerate}
\item For $x \in [0,1]$: \[\sqrt{x} \leq MY(x) \leq \sqrt[3]{x} \]
\item For $x \in [1,+\infty]$: \[\sqrt[3]{x} \leq MY(x) \leq \sqrt{x} \]
\end{enumerate}
\newpage
{\bf Derivative and primitive:}
\begin{enumerate}
\item For $x>0$ the derivative of $M$ is given by the expression:
\[ MY'\left(x\right)=\frac{2}{3 MY^2(x)+2 MY(x)} \] 
\item For $x \geq 0$ a primitive of $MY$ is given by the expression:
\[ \frac{3}{4} x MY\left(x\right)- \frac{x}{12}+ \frac{MY^2\left(x\right)}{24}\]
\end{enumerate}
{\bf Proof of the properties of $MY$}\\
\\
{\bf Equalities:}
\begin{enumerate}
\item 
The objective is to prove:
\[MY(x)=\frac{v}{1-u} \qquad (25)\]
Where
\[ u= 3 \sqrt[3]{\left (\frac {1}{27}-x \right)+\sqrt{x\left (x-\frac {2}{27}\right)}} \qquad v= 3 \sqrt[3]{2x\left(x-\sqrt{x\left (x-\frac {2}{27}\right)}\right)}\]
\[\frac{v}{1-u}=\frac{v(1+u+u^2)}{1-u^3}=\frac{v}{1-u^3}+\frac{vu}{1-u^3}+\frac{vu^2}{1-u^3} \qquad (26)\]
Then notice:
\begin{enumerate}
\item 
\[ \frac{1}{u}= 3\sqrt[3]{\left (\frac {1}{27}-x \right)-\sqrt{x\left (x-\frac {2}{27}\right)}}\]
\item 
\[ v^2= -6 x u\]
\item 
\[ \frac{v^3}{1-u^3}= 2x\]
\end{enumerate}
Therefore, the three terms in $(26)$ can be expressed:
 \[\frac{u v}{1-u^3}=\frac{v\frac{v^2}{(-6x)}}{1-u^3}=-\frac{1}{6x}\frac{v^3}{1-u^3}=-\frac{2x}{6x}=-\frac{1}{3}\]
 \[\frac{u^2 v}{1-u^3}=u\frac{vu}{1-u^3}=-\frac{u}{3}=\sqrt[3]{\left (x-\frac {1}{27} \right)-\sqrt{x\left (x-\frac {2}{27}\right)}}\]
\[\frac{v}{1-u^3}=\frac{1}{u}\frac{vu}{1-u^3}=-\frac{1}{3u}=\sqrt[3]{\left (x-\frac {1}{27} \right)+\sqrt{x\left (x-\frac {2}{27}\right)}}\]
Which adds up to:
\[MY(x)=\frac{v}{1-u}\]
In other words: 
\[MY\left(x\right)= \frac{\sqrt[3]{2x\left(x-\sqrt{x\left (x-\frac {2}{27}\right)}\right)}}{\frac{1}{3}+\sqrt[3]{\left (x-\frac {1}{27} \right)-\sqrt{x\left (x-\frac {2}{27}\right)}}}\]
\item 
For any $x >0 $, define $z=MY(x)$:
\[\frac{z^3+z^2}{2}=x\]
A change of variable $y=1/z$ leads to:
\[y^3-\frac{1}{2x}y-\frac{1}{2x}=0\]
Let's introduce a second change of variable:
 \[y=3 w \sqrt{\frac{1}{6x}}+\sqrt{\frac{1}{6x}}\]
which leads to:
\[\frac{w^3+w^2}{2}=\sqrt{\frac{x}{54}}+\frac{1}{27}\]
so:
 \[ w=MY\left( \sqrt{\frac{x}{54}}+\frac{1}{27} \right)  \]
Since  $y=1/z=\frac{1}{MY(x)}$:
\[MY\left(x\right)\left( 3 MY\left(\sqrt{\frac{x}{54}}+\frac{1}{27}\right)+1\right)=\sqrt{6 x}\]

\item Let $z=MY(x)$:
\[\frac{z^3+z^2}{2}=2x\]
By multiplying both sides by $z^{-5}$:
\[\frac{z^{-3}+z^{-2}}{2}=2x z^{-5}\]
Or
\[\frac{1}{z}=M\left(\frac{x}{z^5}\right)\]
Therefore
\[MY(x)=\frac{1}{MY\left(\frac{x}{(MY(x))^5}\right) } \]

\item For $0 \leq x \leq \frac{2}{27}$, $z_1=MY(x)$ is root of the equation $z^3+z^2=2x$. The other two roots, $z_2$ and $z_3$, can be found by symmetry and using Vieta's formula as provided in section 5.1. Alternatively:
\[ z_2+z_3=-1-z_1 \qquad \textrm{and} \qquad z_2z_3=2x/z_1=z_1(1+z_1)\]
$z_2$ and $z_3$ are therefore roots of a quadratic equation. The discriminant is: 
\[\delta=(1+z_1)^2\frac{1-3 z_1}{1+z_1}\]
Therefore
\[z_2=-(1+z_1)\left(1+\sqrt{\frac{1-3z_1}{1+z_1}}\right)\]
\[z_3=-(1+z_1)\left(1-\sqrt{\frac{1-3z_1}{1+z_1}}\right)\]
Notice $z_2 \leq z_3$. Using the results from section 5.1:\\
\[z_2=-\frac{2}{3}-MY\left(\frac{2}{27}-x\right) =-\frac{1+MY(x)}{2}\left(1+\sqrt{\frac{1-3 MY(x)}{1+MY(x)}}\right)\]
\[z_3=MY\left(\frac{2}{27}-x\right)-MY(x)-\frac{1}{3} =-\frac{1+MY(x)}{2}\left(1-\sqrt{\frac{1-3 MY(x)}{1+MY(x)}}\right)\]
Consequently
\[MY\left(\frac{2}{27}-x\right) =\frac{1+MY(x)}{2}\left(1+\sqrt{\frac{1-3 MY(x)}{1+MY(x)}}\right)-\frac{2}{3}\]
\item $MY$ is the inverse of $f$. Since $\lim \limits_{x \to 0}f(x)=0$ and $\lim \limits_{x \to +\infty} f(x)=+\infty$:
\[\lim \limits_{x \to 0}MY(x)=0 \qquad \textrm{and} \qquad \lim \limits_{x \to +\infty} MY(x)=+\infty\]
Also $MY(x)^2\left(1+MY(x)\right)=2x$ which leads to:
\[\frac{MY(x)}{\sqrt{2x}}=\sqrt{\frac{1}{1+MY(x)}}\]
Therefore:
\[\lim \limits_{x \to 0}\frac{MY(x)}{\sqrt{2x}}=1\]
\\
\\
Similarily, since $MY(x)^3\left(1+\frac{1}{MY(x)}\right)=2x$:
\[\frac{MY(x)}{\sqrt[3]{2x}}=\sqrt[3]{\frac{1}{1+\frac{1}{MY(x)}}}\]
Therefore
\[\lim \limits_{x \to \infty}\frac{MY(x)}{\sqrt[3]{2x}}=1\]
\\
\\
\\
Consequently for any $x > 0$
\[  \sqrt{x}=\lim\limits_{\epsilon \to 0^+} \frac{1}{\epsilon} MY\left(\frac{x \epsilon^2}{2}\right) \qquad \textrm{and} \qquad  \sqrt[3]{x}=\lim\limits_{\epsilon \to 0^+} \epsilon MY\left(\frac{x}{2\epsilon^3}\right)\]
This holds also for $x=0$.
\end{enumerate}
{\bf Inequalities:}
\begin{enumerate}
\item For all positive real numbers $x$
\[f\left(x^{\frac{2}{5}}\right)=\frac{1}{2} \left(x^{\frac{6}{5}}+x^{\frac{4}{5}}\right)=x\left(1+\frac{1}{2} \left(x^{\frac{1}{5}}-x^{\frac{-1}{5}}\right)^2\right) \geq x\]
Since $f$ is strictly increasing
\[MY(x) \leq x^{\frac{2}{5}} \]
In addition
\[MY(x)=\sqrt{\frac{2x}{1+MY(x)}} \geq \sqrt{\frac{2x}{1+x^{\frac{2}{5}}}}\]

\item Let's consider the function $h(x)= x^a$ and define $z=MY(x)$:
\begin{enumerate}
\item For $0 \leq a \leq 1$, $h$ is concave and:
\[\frac{h(z^3)+h(z^2)}{2} \leq h\left(\frac{z^3+z^2}{2}\right)=h(x)\]
Since $h(z^3)=(h(z))^3$ and $h(z^2)=(h(z))^2$:
\[\frac{(h(z))^3+(h(z))^2}{2} \leq h(x) \qquad \textrm{and }\qquad h(z) \leq MY(h(x))\]
Therefore
\[MY\left(x^a\right) \geq \left(MY(x)\right)^a\]
\item Similarily, for $a\leq 0$ or $a\geq 1$, $h$ is convex and we have:
\[MY\left(x^a\right) \leq \left(MY(x)\right)^a\]

\end{enumerate}

\item For all $x \in [0,1]$:
\[f(\sqrt{x}) =x\frac{\sqrt{x}+1}{2} \leq x\]
\[f(\sqrt[3]{x}) =x\frac{x^{-1/3}+1}{2} \geq x\]
Since $f$ is strictly increasing:
\[\sqrt{x} \leq MY(x)\leq \sqrt[3]{x}\]

\item For $x \in [1,+\infty]$
\[f(\sqrt{x}) =x\frac{\sqrt{x}+1}{2} \geq x\]
\[f(\sqrt[3]{x}) =x\frac{x^{-1/3}+1}{2} \leq x\]
Since $f$ is strictly increasing:
\[ \sqrt[3]{x} \leq MY(x) \leq \sqrt{x}\]

\end{enumerate}
{\bf Derivative and primitive:}
\begin{enumerate}
\item Deriving the identity $f(MY(x))=x$ leads to $MY'(x)f'(MY(x))=1$  
\[ MY'\left(x\right)\left(\frac{3}{2}MY^2(x)+MY(x) \right)=1\] 
Which means:
\[ MY'\left(x\right)=\frac{2}{3 MY^2(x)+2 MY(x)} \] 
.
\item Since $MY$ is continuous over $\R^{+}$, a primitive of $MY$ is given by: 
\[\int_{0}^{x} MY(t)dt\]
Let's use a change of variable $MY(t)=w$ (or $t=f(w)$ and $dt=f'(w)dw$):
 \[\int_{0}^{x} MY(t)dt=\int_{0}^{MY(x)} w \left( \frac{3 w^2+2 w}{2}\right)dw\]
Using $x=(MY^3(x)+MY^2(x))/2$, the integral can be simplified to:
\[ \int_{0}^{x} MY(t)dt=\frac{3}{4} x MY\left(x\right)- \frac{x}{12}+ \frac{MY^2\left(x\right)}{24}\]
\end{enumerate}
\section{Acknowledgment}
The authors would like to thank Dr. Nizar Demni, Dr. Hassan Youssfi and Dr. Daniel Moseley for their support and valuable feedback.

\end{document}